\documentclass{article}
\usepackage{amssymb}

\usepackage{graphicx}

\input{tcilatex}

\begin{document}

\begin{center}
\textbf{A self-avoiding curve associated with sums of digits}

F. Oger\bigskip
\end{center}

\bigskip

\noindent \textbf{Abstract.} For each $n\in 
%TCIMACRO{\U{2115} }%
%BeginExpansion
\mathbb{N}
%EndExpansion
^{\ast }$, we write $s_{n}=\left( 1,\ldots ,1,0\right) $\ with $n$ times $1$%
.\ For each $a\in 
%TCIMACRO{\U{2115} }%
%BeginExpansion
\mathbb{N}
%EndExpansion
$, we consider the binary representation $\left( a_{i}\right) _{i\in -%
%TCIMACRO{\U{2115} }%
%BeginExpansion
\mathbb{N}
%EndExpansion
}$\ of $a$ with\ $a_{i}=0$\ for nearly each $i$; we denote by $\alpha _{n}(a)
$ the number of integers $i$ such that $\left( a_{i},\ldots ,a_{i+n}\right)
=s_{n}$. We consider the curve $C_{n}=\left( S_{n,k}\right) _{k\in 
%TCIMACRO{\U{2115} }%
%BeginExpansion
\mathbb{N}
%EndExpansion
^{\ast }}$ which consists of consecutive segments of length $1$ such that,
for each $k\in 
%TCIMACRO{\U{2115} }%
%BeginExpansion
\mathbb{N}
%EndExpansion
^{\ast }$, $S_{n,k+1}$\ is obtained from $S_{n,k}$ by turning right if $%
k+\alpha _{n}(k)-\alpha _{n}(k-1)$\ is even and left otherwise.

$C_{1}$\ is self-avoiding since it is the curve associated to the
alternating folding sequence. In [1], M. Mend\`{e}s France and J. Shallit
conjectured that the curves $C_{n}$ for $n\geq 2$ are also self-avoiding. In
the present paper, we show\ that this property is true for $n=2$. We also
prove that $C_{2}$ has some properties similar to those which were shown in
[2], [3] and [4] for folding curves.\bigskip

For each $a=\left( a_{i}\right) _{i\in -%
%TCIMACRO{\U{2115} }%
%BeginExpansion
\mathbb{N}
%EndExpansion
}\in \left\{ 0,1\right\} ^{-%
%TCIMACRO{\U{2115} }%
%BeginExpansion
\mathbb{N}
%EndExpansion
}$, we define $a+1=\left( b_{i}\right) _{i\in -%
%TCIMACRO{\U{2115} }%
%BeginExpansion
\mathbb{N}
%EndExpansion
}$ as follows:

\noindent - if $a_{i}=1$\ for each $i$, then $b_{i}=0$\ for each $i$;

\noindent - otherwise, there exists $i\in -%
%TCIMACRO{\U{2115} }%
%BeginExpansion
\mathbb{N}
%EndExpansion
$ such that $a_{i}=0$\ and $a_{j}=1$ for $j>i$; we write $b_{j}=a_{j}$ for $%
j<i$, $b_{i}=1$\ and $b_{j}=0$ for $j>i$.

We write $P(a)=a_{0}$. If $a=b+1$ with $b=\left( b_{i}\right) _{i\in -%
%TCIMACRO{\U{2115} }%
%BeginExpansion
\mathbb{N}
%EndExpansion
}$, we denote by $Q\left( a\right) $\ (resp. $R\left( a\right) $) the number
of integers $i$\ such that $\left( b_{i-2},b_{i-1},b_{i}\right) =\left(
1,1,0\right) $ and $\left( a_{i-2},a_{i-1},a_{i}\right) \neq \left(
1,1,0\right) $ (resp. $\left( b_{i-2},b_{i-1},b_{i}\right) \neq \left(
1,1,0\right) $ and $\left( a_{i-2},a_{i-1},a_{i}\right) =\left( 1,1,0\right) 
$). We have $P(a),Q(a),R(a)\in \left\{ 0,1\right\} $.

Each positive integer is represented by such an $a$ with $a_{i}=1$\ for
finitely many integers $i$.

We write $e=\left( 1,0\right) $ and $f=\left( 0,1\right) $. We define two
types of curves consisting of consecutive segments of length $1$:

For $a=\left( a_{i}\right) _{i\in -%
%TCIMACRO{\U{2115} }%
%BeginExpansion
\mathbb{N}
%EndExpansion
}\in \left\{ 0,1\right\} ^{-%
%TCIMACRO{\U{2115} }%
%BeginExpansion
\mathbb{N}
%EndExpansion
}$ with $\left\{ i\in -%
%TCIMACRO{\U{2115} }%
%BeginExpansion
\mathbb{N}
%EndExpansion
\mid a_{i}=0\right\} $ and $\left\{ i\in -%
%TCIMACRO{\U{2115} }%
%BeginExpansion
\mathbb{N}
%EndExpansion
\mid a_{i}=1\right\} $ infinite, for $x\in 
%TCIMACRO{\U{2124} }%
%BeginExpansion
\mathbb{Z}
%EndExpansion
^{2}$ and for $g\in \left\{ e,-e,f,-f\right\} $, we write $C_{a,x,g}=\left(
X_{k}\right) _{k\in 
%TCIMACRO{\U{2124} }%
%BeginExpansion
\mathbb{Z}
%EndExpansion
}$ with $X_{1}=\left[ x,x+g\right] $ and, for each $k\in 
%TCIMACRO{\U{2124} }%
%BeginExpansion
\mathbb{Z}
%EndExpansion
$, $X_{k+1}$ obtained from $X_{k}$ by turning right if $P(a+k)+Q(a+k)+R(a+k)$
is even, and left otherwise.

For $g=\mp e$ and $h=\mp f$, or $g=\mp f$ and $h=\mp e$, we write $%
C_{g,h}=\left( X_{k}\right) _{k\in 
%TCIMACRO{\U{2124} }%
%BeginExpansion
\mathbb{Z}
%EndExpansion
}$ with $X_{0}=\left[ -h,0\right] $, $X_{1}=\left[ 0,g\right] $ and, for
each $k\in 
%TCIMACRO{\U{2124} }%
%BeginExpansion
\mathbb{Z}
%EndExpansion
^{\ast }$, $X_{k+1}$ obtained from $X_{k}$ by turning right if $%
P(k)+Q(k)+R(k)$ is even, and left otherwise. We write $C=C_{f,e}$.

For any $X,Y\subset 
%TCIMACRO{\U{211d} }%
%BeginExpansion
\mathbb{R}
%EndExpansion
^{2}$, an \emph{isomorphism} from $X$ to $Y$ is a translation $\tau $ such
that $\tau (X)=Y$.\ They are \emph{locally isomorphic} if, for each $x\in 
%TCIMACRO{\U{211d} }%
%BeginExpansion
\mathbb{R}
%EndExpansion
^{2}$\ (resp. $y\in 
%TCIMACRO{\U{211d} }%
%BeginExpansion
\mathbb{R}
%EndExpansion
^{2}$) and each $s\in 
%TCIMACRO{\U{211d} }%
%BeginExpansion
\mathbb{R}
%EndExpansion
^{+}$, there exists $y\in 
%TCIMACRO{\U{211d} }%
%BeginExpansion
\mathbb{R}
%EndExpansion
^{2}$\ (resp. $x\in 
%TCIMACRO{\U{211d} }%
%BeginExpansion
\mathbb{R}
%EndExpansion
^{2}$) such that $(B(x,s)\cap X,x)\cong (B(y,s)\cap Y,y)$.

We say that $X\subset 
%TCIMACRO{\U{211d} }%
%BeginExpansion
\mathbb{R}
%EndExpansion
^{2}$ satisfies the \emph{local isomorphism property} if, for any $x\in 
%TCIMACRO{\U{211d} }%
%BeginExpansion
\mathbb{R}
%EndExpansion
^{2}$ and $s\in 
%TCIMACRO{\U{211d} }%
%BeginExpansion
\mathbb{R}
%EndExpansion
^{+}$, there exists $t\in 
%TCIMACRO{\U{211d} }%
%BeginExpansion
\mathbb{R}
%EndExpansion
^{+}$ such that each $B(y,t)$ with $y\in 
%TCIMACRO{\U{211d} }%
%BeginExpansion
\mathbb{R}
%EndExpansion
^{2}$\ contains some $z$ with $(B(z,s)\cap X,z)\cong (B(x,s)\cap X,x)$.

We say that a curve $D\subset 
%TCIMACRO{\U{211d} }%
%BeginExpansion
\mathbb{R}
%EndExpansion
^{2}$ satisfies the \emph{weak local isomorphism property} if, for any $x\in
D$ and $s\in 
%TCIMACRO{\U{211d} }%
%BeginExpansion
\mathbb{R}
%EndExpansion
^{+}$, there exists $t\in 
%TCIMACRO{\U{211d} }%
%BeginExpansion
\mathbb{R}
%EndExpansion
^{+}$ such that each $B(y,t)$ with $y\in D$\ contains some $z$ with $%
(B(z,s)\cap D,z)\cong (B(x,s)\cap D,x)$. It follows that a point $x\in D$
cannot be distinguished by the properties of $D$ around $x$.

We are going to prove the following results:\bigskip

\noindent \textbf{Theorem 1}. Each curve $C_{a,x,g}$\ or $C_{g,h}$ is
self-avoiding and satisfies the weak local isomorphism property.\bigskip

We note that such a curve cannot satisfy the local isomorphism property
since, for each $s\in 
%TCIMACRO{\U{211d} }%
%BeginExpansion
\mathbb{R}
%EndExpansion
^{+}$, there exists $x\in 
%TCIMACRO{\U{211d} }%
%BeginExpansion
\mathbb{R}
%EndExpansion
^{2}$ such that $B(x,s)$ contains no point of that curve.\bigskip

\noindent \textbf{Theorem 2}. The set consisting of the curves $C_{a,x,g}$\
and the curves $C_{g,h}$ is the union of two local isomorphism classes. Each
curve which is locally isomorphic to one of them is isomorphic to another
one. We pass from one class to the other one by replacing $\mp e$ with $\mp
f $ and $\mp f$ with $\mp e$.\bigskip

We consider the monoid $\Omega $\ generated by $u,\overline{u},v,\overline{v}%
,w,\overline{w}$ and the endomorphism $\varphi $\ of $\Omega $\ such that $%
\varphi \left( u\right) =uv$, $\varphi \left( \overline{u}\right) =\overline{%
u}\overline{v}$, $\varphi \left( v\right) =uw$, $\varphi \left( \overline{v}%
\right) =\overline{u}\overline{w}$, $\varphi \left( w\right) =\overline{u}w$%
, $\varphi \left( \overline{w}\right) =u\overline{w}$.

We have $\varphi ^{2}(u)=uvuw$ and $\varphi ^{3}(u)=uvuwuv\overline{u}w$.
For each $n\in 
%TCIMACRO{\U{2115} }%
%BeginExpansion
\mathbb{N}
%EndExpansion
$, $\varphi ^{n}(u)$ is an initial segment of $\varphi ^{n+1}(u)$.

$\varphi $\ commutes with the endomorphism $x\rightarrow \overline{x}$\
where $\overline{\overline{u}}=u$, $\overline{\overline{v}}=v$ and $%
\overline{\overline{w}}=w$.

We consider the sequences $\left( u_{n}\right) _{n\in 
%TCIMACRO{\U{2115} }%
%BeginExpansion
\mathbb{N}
%EndExpansion
},\left( v_{n}\right) _{n\in 
%TCIMACRO{\U{2115} }%
%BeginExpansion
\mathbb{N}
%EndExpansion
},\left( w_{n}\right) _{n\in 
%TCIMACRO{\U{2115} }%
%BeginExpansion
\mathbb{N}
%EndExpansion
}\subset 
%TCIMACRO{\U{2124} }%
%BeginExpansion
\mathbb{Z}
%EndExpansion
^{2}$ with $u_{0}=e$, $v_{0}=f$, $w_{0}=f$ and, for each $n\in 
%TCIMACRO{\U{2115} }%
%BeginExpansion
\mathbb{N}
%EndExpansion
$, $u_{n+1}=u_{n}+v_{n}$, $v_{n+1}=u_{n}+w_{n}$, $w_{n+1}=w_{n}-u_{n}$. We
have

\noindent $u_{1}=e+f$, $v_{1}=e+f$, $w_{1}=-e+f$,

\noindent $u_{2}=2e+2f$, $v_{2}=2f$, $w_{2}=-2e$,

\noindent $u_{3}=2e+4f$, $v_{3}=2f$, $w_{3}=-4e-2f$,

\noindent $u_{4}=2e+6f$, $v_{4}=-2e+2f$, $w_{4}=-6e-6f$.

We see by induction on $n$ that, for $n\geq 2$, $u_{n}$, $v_{n}$ and $w_{n}$
are pairwise non colinear and, by rotating around $0$ in anticlockwise
direction, we successively find $u_{n}$, $v_{n}$, $w_{n}$, $-u_{n}$, $-v_{n}$%
, $-w_{n}$.

We consider the homomorphism $\psi $\ from $\Omega $\ to $\left( 
%TCIMACRO{\U{2124} }%
%BeginExpansion
\mathbb{Z}
%EndExpansion
^{2},+\right) $\ such that $\psi (u)=e$, $\psi (\overline{u})=-e$, $\psi
(v)=f$, $\psi (\overline{v})=-f$, $\psi (w)=f$, $\psi (\overline{w})=-f$. It
follows from the definition of $\varphi $ that we have $\psi \left( \varphi
^{n}\left( u\right) \right) =u_{n}$, $\psi \left( \varphi ^{n}\left(
v\right) \right) =v_{n}$ and $\psi \left( \varphi ^{n}\left( w\right)
\right) =w_{n}$ for each $n\in 
%TCIMACRO{\U{2115} }%
%BeginExpansion
\mathbb{N}
%EndExpansion
$.

For each $n\in 
%TCIMACRO{\U{2115} }%
%BeginExpansion
\mathbb{N}
%EndExpansion
$, we consider the curve $C_{n}=\left( S_{1},\ldots ,S_{2^{n}}\right) $
starting at $0$\ which consists of $2^{n}$ consecutive segments of length $1$
defined as follows: we write $\varphi ^{n}(u)=x_{1}\cdots x_{2^{n}}$ with $%
x_{1},\ldots ,x_{2^{n}}\in \left\{ u,\overline{u},v,\overline{v},w,\overline{%
w}\right\} $; for $1\leq i\leq 2^{n}$, we take $S_{i}$ isomorphic to $u_{0}$
(resp.\ $-u_{0}$, $v_{0}$,\ $-v_{0}$, $w_{0}$, $-w_{0}$) if $x_{i}=u$\
(resp. $\overline{u}$, $v$, $\overline{v}$, $w$, $\overline{w}$). For each $%
n\in 
%TCIMACRO{\U{2115} }%
%BeginExpansion
\mathbb{N}
%EndExpansion
$, $C_{n}$ is an initial segment of $C_{n+1}$. The curve $C_{9}$\ is shown
in Figure 1 below.\bigskip

\begin{center}
\includegraphics[scale=0.75]{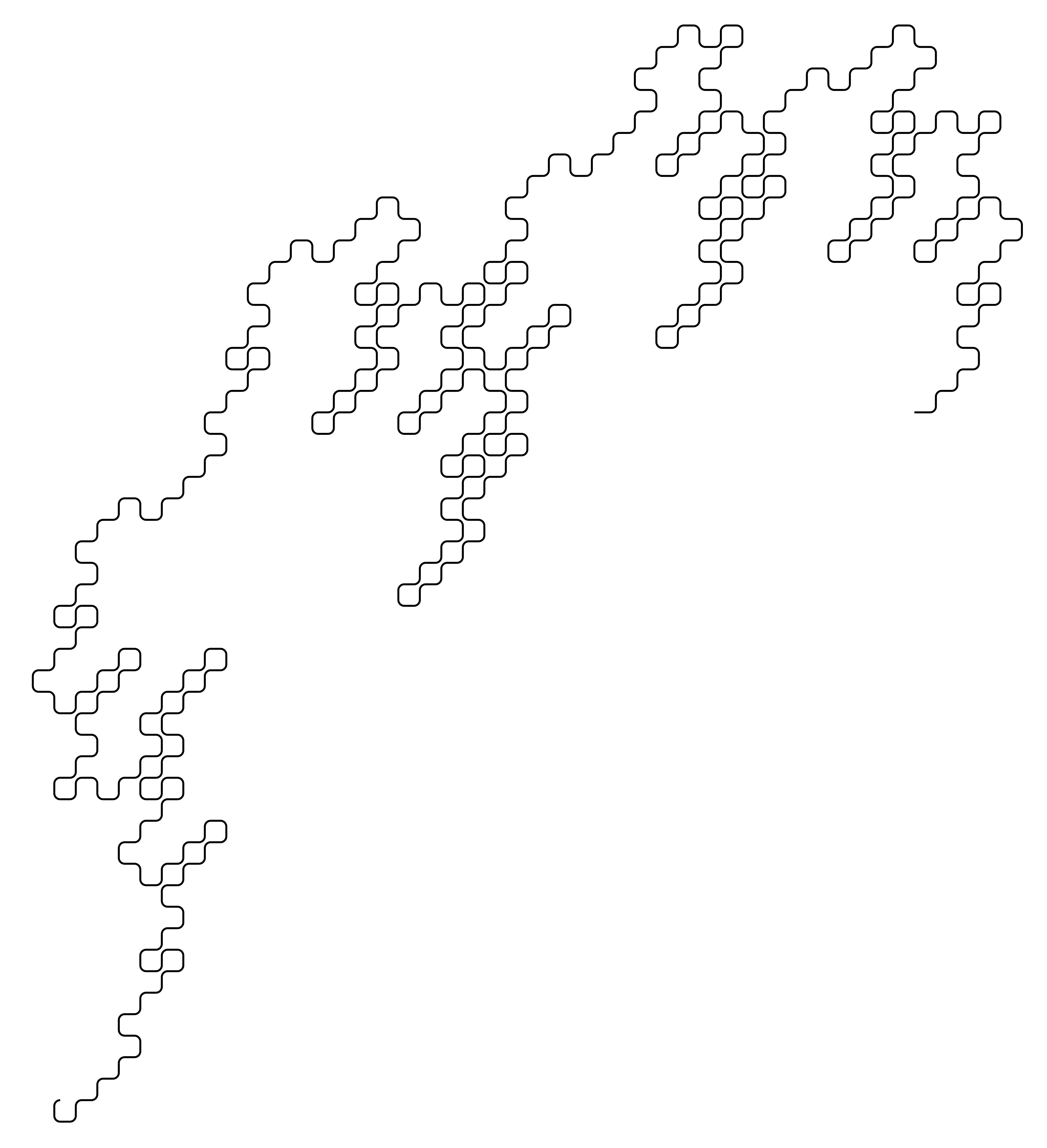}

Figure 1\bigskip
\end{center}

\noindent \textbf{Lemma 3}. For each $n\in 
%TCIMACRO{\U{2115} }%
%BeginExpansion
\mathbb{N}
%EndExpansion
$, $C_{n}$ is an initial segment of the part of $C$ starting at $0$.\bigskip

\noindent \textbf{Proof.} We prove this result by induction on $n$. We see
from Figure 1 that it is true for $n\leq 4$. Now we show that, if it is true
for $n\geq 4$, then it is also true\ for $n+1$.

We write $C_{n+1}=(X_{i})_{1\leq i\leq 2^{n+1}}$\ and $C=(Y_{i})_{i\in 
%TCIMACRO{\U{2124} }%
%BeginExpansion
\mathbb{Z}
%EndExpansion
}$. By the induction hypothesis, $X_{i}=Y_{i}$\ \ is true for\ $1\leq i\leq
2^{n}$.\ We must prove that it is also true\ for\ $2^{n}+1\leq i\leq 2^{n+1}$%
.

The equalities $\varphi ^{3}(u)=uvuwuv\overline{u}w=\varphi (u)uw\varphi (u)%
\overline{u}w$ and $\overline{\varphi ^{n-2}(u)}=\varphi ^{n-2}(\overline{u}%
) $ imply $\varphi ^{n+1}(u)=\varphi ^{n-1}(u)\varphi ^{n-2}(u)\varphi
^{n-2}(w)\varphi ^{n-1}(u)\overline{\varphi ^{n-2}(u)}\varphi ^{n-2}(w)$.

For $k\geq 2$, the words $\varphi ^{k}(u)$, $\overline{\varphi ^{k}(u)}$, $%
\varphi ^{k}(v)$, $\varphi ^{k}(w)$ respectively begin with $u$, $\overline{u%
}$, $u$, $\overline{u}$\ and finish with $w$, $\overline{w}$, $w$, $w$. It
follows that we turn right from $X_{2^{n}}$\ to $X_{2^{n}+1}$, right from $%
X_{2^{n-1}}$\ to $X_{2^{n-1}+1}$, left from $X_{2^{n}+2^{n-1}}$\ to $%
X_{2^{n}+2^{n-1}+1}$, left from $X_{2^{n-1}+2^{n-2}}$\ to $%
X_{2^{n-1}+2^{n-2}+1}$\ and right from $X_{2^{n}+2^{n-1}+2^{n-2}}$\ to $%
X_{2^{n}+2^{n-1}+2^{n-2}+1}$. On the other hand, for each $k\in \left\{
1,\ldots ,2^{n}-1\right\} -\left\{ 2^{n-1},2^{n-1}+2^{n-2}\right\} $, if we
turn right (resp. left) from $X_{k}$\ to $X_{k+1}$, then we turn right
(resp. left) from $X_{2^{n}+k}$\ to $X_{2^{n}+k+1}$.

It follows from the definition of $C$ that these properties are also true
for $(Y_{i})_{1\leq i\leq 2^{n+1}}$. Actually,\ for $0\leq k\leq 2^{n-1}-1$\
or $2^{n-1}+2^{n-2}\leq k\leq 2^{n}-1$, the representations of $k$ and $%
2^{n}+k$ in binary notation contain the same number of sequences $110$. For $%
2^{n-1}\leq k\leq 2^{n-1}+2^{n-2}-1$, the representation of $2^{n}+k$
contains one more sequence $110$ than\ the representation of $k$. The
representations of $2^{n}-1$ and $2^{n}$ contain no sequence $110$.

It follows that we also have $X_{i}=Y_{i}$\ for\ $2^{n}+1\leq i\leq 2^{n+1}$%
.~~$\blacksquare $\bigskip

For any integers $2\leq m\leq n$, we denote by $C_{m,n}$\ the curve\ with $%
2^{n-m}$ segments defined\ as follows: we group the $2^{n}$ segments of $%
C_{n}$\ in $2^{n-m}$ sequences of $2^{m}$ consecutive segments and we
replace each of these sequences with the segment which joins its endpoints.

For $2\leq m\leq n$, we write $\varphi ^{n}(u)=s_{1}\cdots s_{2^{n-m}}$
with\ $s_{1},\ldots ,s_{2^{n-m}}\in \left\{ u,\overline{u},v,\overline{v},w,%
\overline{w}\right\} ^{2^{m}}$ and $C_{m,n}=\left( S_{1},\ldots
,S_{2^{n-m}}\right) $ with\ $S_{1},\ldots ,S_{2^{n-m}}$\ consecutive
segments. For $1\leq i\leq 2^{n}$, we have $\psi \left( s_{i}\right)
=y_{i}-x_{i}$ where $S_{i}=\left[ x_{i},y_{i}\right] $.\ In particular, each
segment of $C_{m,n}$\ is isomorphic to one of the vectors $u_{m}$,\ $-u_{m}$%
, $v_{m}$,\ $-v_{m}$, $w_{m}$, $-w_{m}$.

For $x,y,z\in 
%TCIMACRO{\U{211d} }%
%BeginExpansion
\mathbb{R}
%EndExpansion
^{2}$ with $y,z$\ not colinear, we denote by $P_{x,y,z}$ the parallelogram
with vertices $x$, $x+y$, $x+z$, $x+y+z$.\bigskip

\noindent \textbf{Lemma 4}. For $2\leq m\leq n$, $C_{m,n}$ is self-avoiding
and there exists a tiling $\Pi _{m,n}$ of $%
%TCIMACRO{\U{211d} }%
%BeginExpansion
\mathbb{R}
%EndExpansion
^{2}$ by parallelograms $P_{x,y,z}$ with $x\in 
%TCIMACRO{\U{2124} }%
%BeginExpansion
\mathbb{Z}
%EndExpansion
^{2}$\ and $y,z\in \left\{ u_{m},v_{m},w_{m}\right\} $\ such that:

\noindent 1) there exists no $x\in 
%TCIMACRO{\U{2124} }%
%BeginExpansion
\mathbb{Z}
%EndExpansion
^{2}$\ such that $P_{x,u_{m},v_{m}},P_{x+v_{m},u_{m},v_{m}}\in \Pi _{m,n}$
or\ such that $P_{x,v_{m},w_{m}},P_{x+v_{m},v_{m},w_{m}}\in \Pi _{m,n}$;

\noindent 2) for each segment $S$ of $C_{m,n}$:

\noindent if $S\cong u_{m}$, then there exists\ $P_{x,u_{m},w_{m}}\in \Pi
_{m,n}$\ such that $S=\left[ x+w_{m},x+w_{m}+u_{m}\right] $;

\noindent if $S\cong -u_{m}$, then there exists\ $P_{x,u_{m},w_{m}}\in \Pi
_{m,n}$\ such that $S=\left[ x+u_{m},x\right] $;

\noindent if $S\cong v_{m}$, then there exists\ $P_{x,v_{m},w_{m}}\in \Pi
_{m,n}$\ such that $S=\left[ x+w_{m},x+w_{m}+v_{m}\right] $;

\noindent if $S\cong -v_{m}$, then there exists\ $P_{x,v_{m},w_{m}}\in \Pi
_{m,n}$\ such that $S=\left[ x+v_{m},x\right] $;

\noindent if $S\cong w_{m}$, then there exists\ $P_{x,u_{m},w_{m}}\in \Pi
_{m,n}$\ such that $S=\left[ x+u_{m},x+u_{m}+w_{m}\right] $, or $%
P_{x,v_{m},w_{m}}\in \Pi _{m,n}$\ such that $S=\left[ x+v_{m},x+v_{m}+w_{m}%
\right] $;

\noindent if $S\cong -w_{m}$, then there exists\ $P_{x,u_{m},w_{m}}\in \Pi
_{m,n}$\ such that $S=\left[ x+w_{m},x\right] $, or $P_{x,v_{m},w_{m}}\in
\Pi _{m,n}$\ such that $S=\left[ x+w_{m},x\right] $.\bigskip

\noindent \textbf{Proof.} $C_{n,n}$ is self-avoiding since it just consists
of the\ segment $\left[ 0,u_{n}\right] $. We take for $\Pi _{n,n}$ the
regular tiling by parallelograms $P_{x,u_{n},w_{n}}$ such that $\left[
0,u_{n}\right] $ is an edge of one of them.

Now we show that, for $2\leq m\leq n-1$, if Lemma 4 is true for $m+1$, then
it is also true for $m$. Figures 2, 3, 4 below are exact for $m=3$. For
other values of $m$, the parallelograms have different shapes and dimensions.

First we prove that $C_{m,n}$ is self-avoiding.

$C_{m,n}$ is obtained from $C_{m+1,n}$ by replacing each segment $S$ with a
pair of consecutive segments $\left( S_{1},S_{2}\right) $\ such that:

\noindent if $S\cong u_{m+1}$, then $S_{1}\cong u_{m}$ and $S_{2}\cong v_{m}$%
;

\noindent if $S\cong -u_{m+1}$, then $S_{1}\cong -u_{m}$ and $S_{2}\cong
-v_{m}$;

\noindent if $S\cong v_{m+1}$, then $S_{1}\cong u_{m}$ and $S_{2}\cong w_{m}$%
;

\noindent if $S\cong -v_{m+1}$, then $S_{1}\cong -u_{m}$ and $S_{2}\cong
-w_{m}$;

\noindent if $S\cong w_{m+1}$, then $S_{1}\cong -u_{m}$ and $S_{2}\cong
w_{m} $;

\noindent if $S\cong -w_{m+1}$, then $S_{1}\cong u_{m}$ and $S_{2}\cong
-w_{m}$.

\noindent Moreover, as it is shown in Figure 2, Lemma 4 for $m+1$ implies
that, for any such $S,S_{1},S_{2}$, there exists $P\in \Pi _{m+1,n}$ such
that $S$ is a side of $P$ and such that $S_{1},S_{2}$ are in the interior of 
$P$.

Now, consider two distinct pairs of consecutive segments $\left(
S_{1},S_{2}\right) ,\left( T_{1},T_{2}\right) \subset C_{m,n}$, obtained
from two segments $S,T\in C_{m+1,n}$, and contained in two parallelograms $%
P,Q\in \Pi _{m+1,n}$.

First suppose $P=Q$. As $C_{m+1,n}$ is self-avoiding, the orientations of $S$
and $T$ define the same direction of rotation around the center of $P$. It
follows that $S,T$ are on opposite sides of $P$, and $\left(
S_{1},S_{2}\right) $, $\left( T_{1},T_{2}\right) $\ have no common point, or
they have just one common point, which is the terminal point of $S_{1},T_{1}$
and the initial point of $S_{2},T_{2}$, and they do not cross each other at
that point.

Now suppose $P\neq Q$. Then $\left( S_{1},S_{2}\right) $, $\left(
T_{1},T_{2}\right) $\ have no common point, except if the initial point of $%
S $ is the initial or the terminal point of $T$, or if the terminal point of 
$S $ is the initial or the terminal point of $T$. In that case, $\left(
S_{1},S_{2}\right) $, $\left( T_{1},T_{2}\right) $\ only have that common
point since $S\neq T$.

It remains to be proved that, for each $x\in 
%TCIMACRO{\U{2124} }%
%BeginExpansion
\mathbb{Z}
%EndExpansion
^{2}$\ and any distinct pairs of consecutive segments $\left(
S_{1},S_{2}\right) ,\left( S_{3},S_{4}\right) \subset C_{m,n}$, with $x$\
terminal point of $S_{1},S_{3}$ and initial point of $S_{2},S_{4}$, if $x$
is a vertex of the parallelograms belonging to $\Pi _{m+1,n}$ and containing 
$S_{1},S_{2},S_{3},S_{4}$, then $S_{1},S_{2},S_{3},S_{4}$\ are distinct and
the curves $\left( S_{1},S_{2}\right) ,\left( S_{3},S_{4}\right) $ do not
cross each other in $x$.

We have $\left( T_{1},T_{2}\right) \neq \left( T_{3},T_{4}\right) $ for the
segments $T_{1},T_{2},T_{3},T_{4}\in C_{m+1,n}$\ from which $%
S_{1},S_{2},S_{3},S_{4}$ are obtained. As $C_{m+1,n}$\ is self-avoiding, $%
T_{1},T_{2},T_{3},T_{4}$\ are distinct\ and $\left( T_{1},T_{2}\right)
,\left( T_{3},T_{4}\right) $ do not cross each other in $x$. The
parallelograms $P_{1},P_{2},P_{3},P_{4}\in \Pi _{m+1,n}$\ associated to $%
T_{1},T_{2},T_{3},T_{4}$\ are also distinct\ and $x$ is their common vertex.
As $S_{1},S_{2},S_{3},S_{4}$ are respectively in the interiors of $%
P_{1},P_{2},P_{3},P_{4}$, it follows that $S_{1},S_{2},S_{3},S_{4}$\ are
distinct and $\left( S_{1},S_{2}\right) $, $\left( S_{3},S_{4}\right) $ do
not cross each other in $x$.

We define $\Pi _{m,n}$ from $\Pi _{m+1,n}$ (see Figure 3) by replacing:

\noindent -\ each $P_{x,u_{m+1},v_{m+1}}$ with the segments $\left[ x,x+u_{m}%
\right] $, $\left[ x+u_{m},x+u_{m}+v_{m}\right] $,

\noindent $\left[ x+u_{m},x+u_{m}+w_{m}\right] $, $\left[
x+v_{m+1},x+v_{m+1}+v_{m}\right] $,

\noindent $\left[ x+v_{m+1}+v_{m},x+v_{m+1}+v_{m}+u_{m}\right] $, $\left[
x+v_{m+1}+v_{m},x+v_{m+1}+v_{m}-w_{m}\right] $;

\noindent - each\ $P_{x,u_{m+1},w_{m+1}}$ with the segments $\left[ x,x+w_{m}%
\right] $, $\left[ x+w_{m},x+w_{m}-u_{m}\right] $, $\left[
x+w_{m},x+w_{m}+v_{m}\right] $, $\left[ x,x+v_{m}\right] $, $\left[
x+v_{m},x+v_{m}+u_{m}\right] $, $\left[ x+v_{m},x+v_{m}+w_{m}\right] $;

\noindent -\ each $P_{x,v_{m+1},w_{m+1}}$ with the segments $\left[ x,x+w_{m}%
\right] $, $\left[ x+w_{m},x+w_{m}-u_{m}\right] $,

\noindent $\left[ x+w_{m},x+w_{m}+u_{m}\right] $, $\left[
x+w_{m},x+w_{m}+w_{m}\right] $.

Each tile of $\Pi _{m,n}$ is obtained from $1$ or $2$ tiles of $\Pi _{m+1,n}$%
. Figure 4 shows, for each possible position of two adjacent tiles $P,Q\in
\Pi _{m+1,n}$, the tiles of $\Pi _{m,n}$ which are obtained from $P$ and $Q$%
. We see from Figure 4 that $\Pi _{m,n}$ is a tiling of $%
%TCIMACRO{\U{211d} }%
%BeginExpansion
\mathbb{R}
%EndExpansion
^{2}$ by parallelograms $P_{x,y,z}$ with $x\in 
%TCIMACRO{\U{2124} }%
%BeginExpansion
\mathbb{Z}
%EndExpansion
^{2}$\ and $y,z\in \left\{ u_{m},v_{m},w_{m}\right\} $, and that there
exists no $x\in 
%TCIMACRO{\U{2124} }%
%BeginExpansion
\mathbb{Z}
%EndExpansion
^{2}$\ such that $P_{x,u_{m},v_{m}},P_{x+v_{m},u_{m},v_{m}}\in \Pi _{m,n}$,
or\ such that $P_{x,v_{m},w_{m}},P_{x+v_{m},v_{m},w_{m}}\in \Pi _{m,n}$.

We see from Figure 2 that, for each segment $S$ of $C_{m,n}$:

\noindent if $S\cong u_{m}$, then there exists\ $P_{x,u_{m},w_{m}}\in \Pi
_{m,n}$\ such that $S=\left[ x+w_{m},x+w_{m}+u_{m}\right] $;

\noindent if $S\cong -u_{m}$, then there exists\ $P_{x,u_{m},w_{m}}\in \Pi
_{m,n}$\ such that $S=\left[ x+u_{m},x\right] $;

\noindent if $S\cong v_{m}$, then there exists\ $P_{x,v_{m},w_{m}}\in \Pi
_{m,n}$\ such that $S=\left[ x+w_{m},x+w_{m}+v_{m}\right] $;

\noindent if $S\cong -v_{m}$, then there exists\ $P_{x,v_{m},w_{m}}\in \Pi
_{m,n}$\ such that $S=\left[ x+v_{m},x\right] $;

\noindent if $S\cong w_{m}$, then there exists\ $P_{x,u_{m},w_{m}}\in \Pi
_{m,n}$\ such that $S=\left[ x+u_{m},x+u_{m}+w_{m}\right] $, or $%
P_{x,v_{m},w_{m}}\in \Pi _{m,n}$\ such that $S=\left[ x+v_{m},x+v_{m}+w_{m}%
\right] $;

\noindent if $S\cong -w_{m}$, then there exists\ $P_{x,u_{m},w_{m}}\in \Pi
_{m,n}$\ such that $S=\left[ x+w_{m},x\right] $, or $P_{x,v_{m},w_{m}}\in
\Pi _{m,n}$\ such that $S=\left[ x+w_{m},x\right] $.~~$\blacksquare $\bigskip

\begin{center}
\includegraphics[scale=0.6]{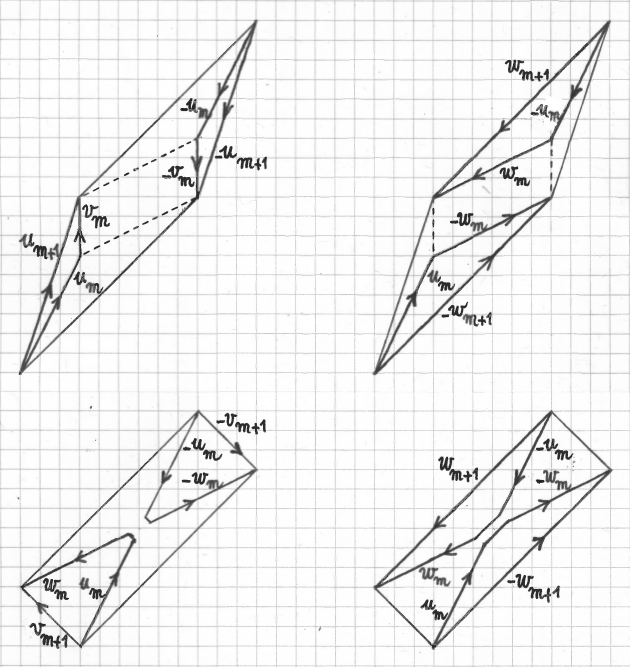}

Figure 2\bigskip

\includegraphics[scale=1.05]{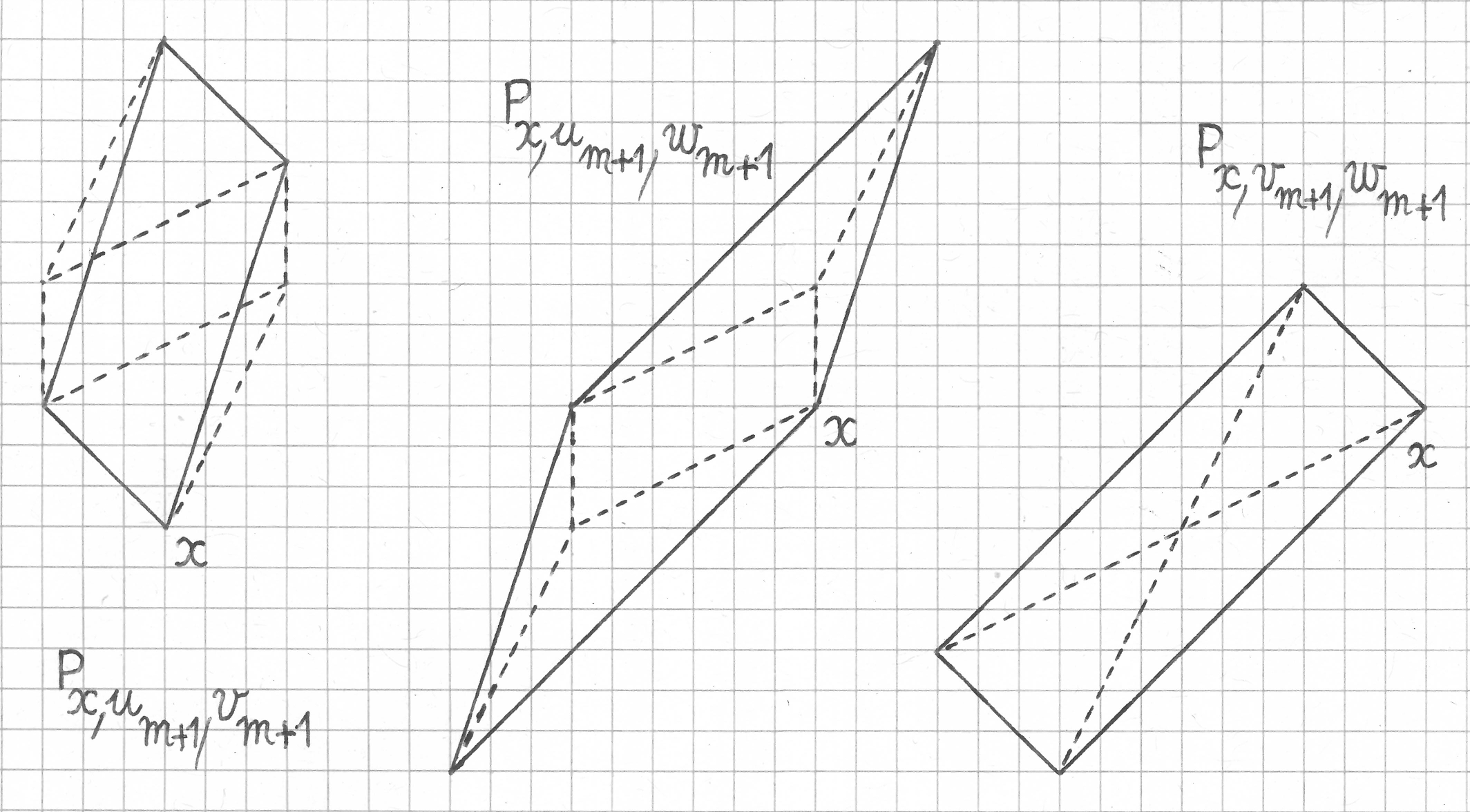}

Figure 3\bigskip

\includegraphics[scale=0.42]{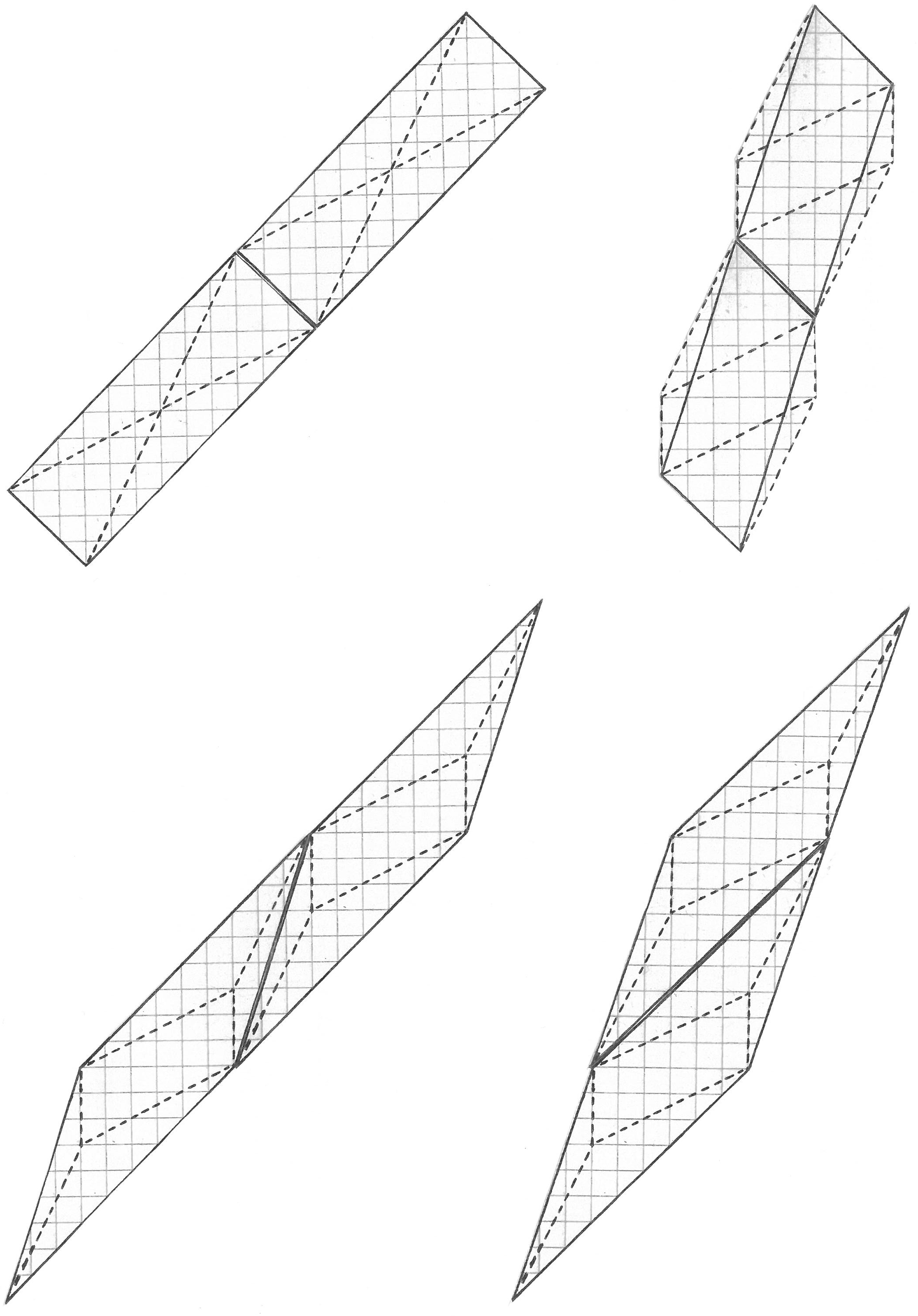}

\includegraphics[scale=0.42]{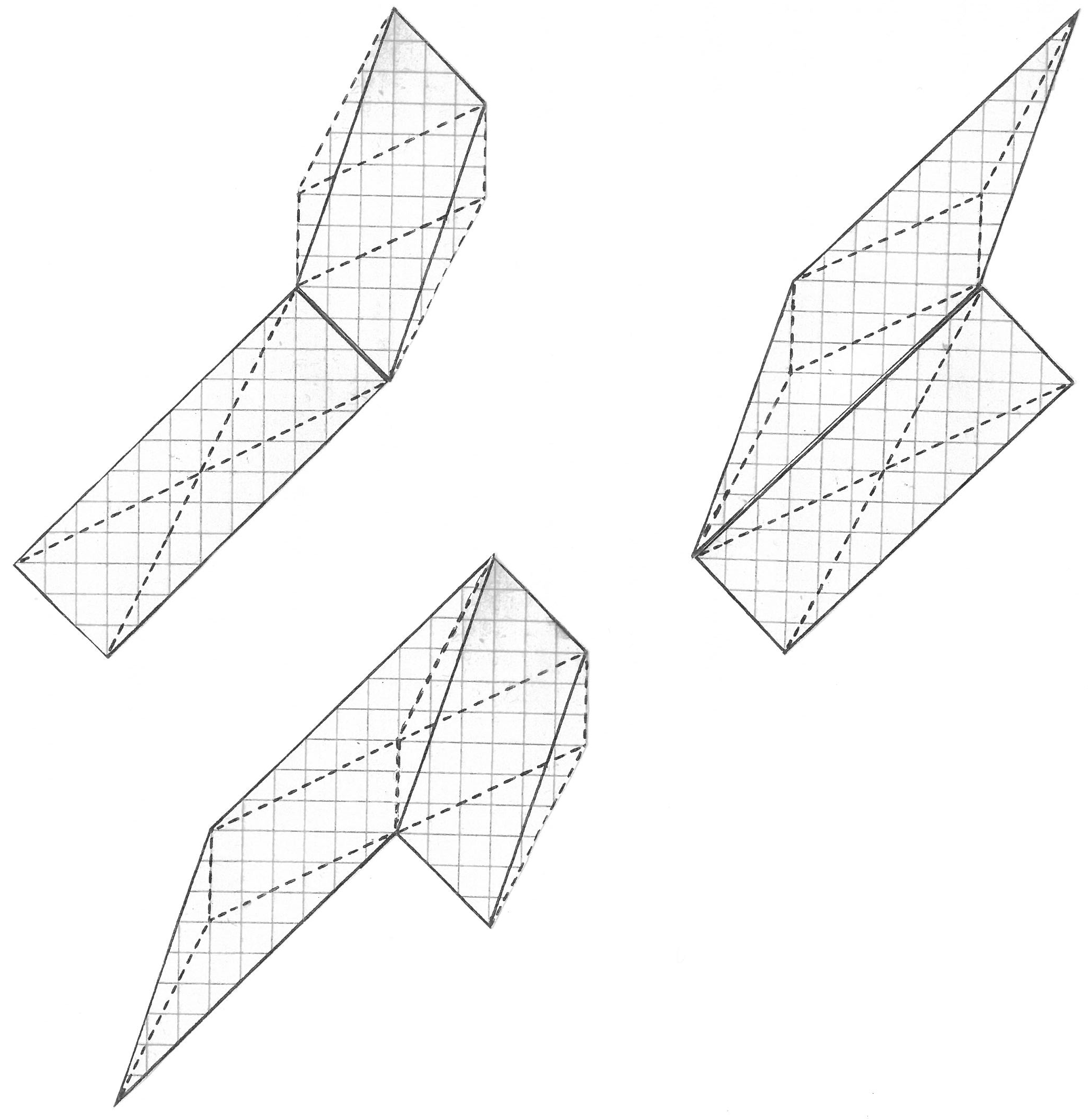}

Figure 4\bigskip
\end{center}

\noindent \textbf{Lemma 5}. Each $C_{n}$\ is self-avoiding.\bigskip

\noindent \textbf{Proof.} In the proof of\ Lemma 4, we first showed that,
for $2\leq m\leq n-1$, Lemma 4 for $C_{m+1,n}$ implies $C_{m,n}$
self-avoiding. Here we show that, for $n\geq 2$, Lemma 4 for $C_{2,n}$
implies $C_{n}$ self-avoiding. The proof is exactly the same, except that
each segment of $C_{2,n}$ is replaced with a sequence of $4$ consecutive
segments, while each segment of $C_{m+1,n}$ was replaced with a pair of
consecutive segments (see Figure 5 below).~~$\blacksquare $\bigskip

\begin{center}
\includegraphics[scale=1.2]{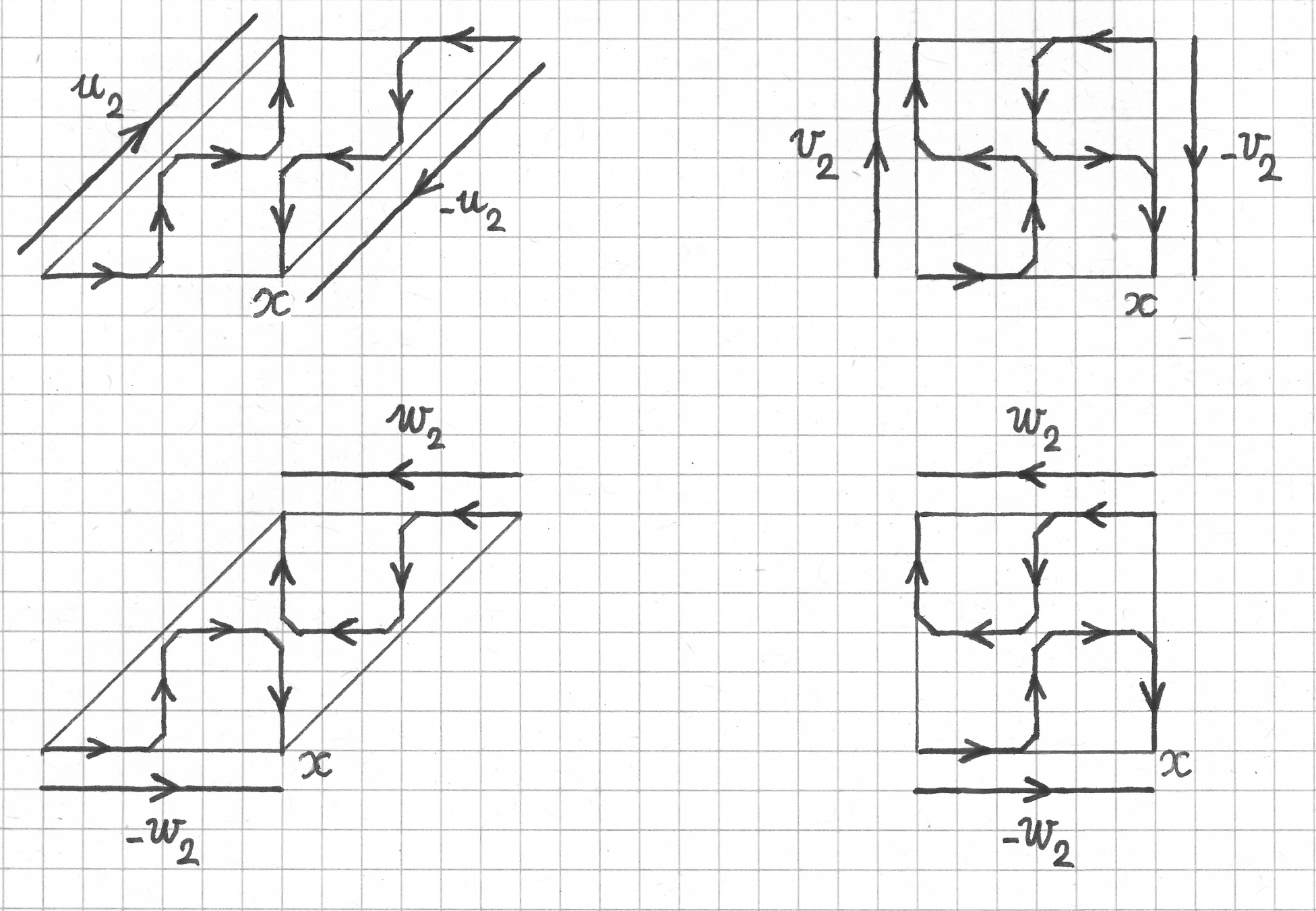}

Figure 5\bigskip
\end{center}

\noindent \textbf{Proof of Theorems 1 and 2.} First we consider two curves $%
A=\left( X_{k}\right) _{k\in 
%TCIMACRO{\U{2124} }%
%BeginExpansion
\mathbb{Z}
%EndExpansion
}$ and $B=\left( Y_{k}\right) _{k\in 
%TCIMACRO{\U{2124} }%
%BeginExpansion
\mathbb{Z}
%EndExpansion
}$, each of them equal to some $C_{g,h}$ or some $C_{a,x,g}$ with $a=(\ldots
,\alpha _{i},\ldots ,\alpha _{1},0)$. We show that, if $X_{k}$\ and $Y_{k}$\
are isomorphic to $\mp e$\ for $k$ odd,\ then each subcurve of $A$ with
finitely many segments appears everywhere in $B$.

This condition can always be realized by rotating $A$, or $B$, or $A$ and $B$%
, with a $\pi /2$ angle, if necessary. If it is only realized for one of the
curves $A,B$, then $A$ and $B$ are not locally isomorphic since only one of
them contains some consecutive segments with endpoints $y$, $y+e$, $y+e+f$, $%
y+2e+f$, $y+2e+2f$, $y+3e+2f$, $y+3e+3f$.

There exist some sequences $a_{k}=\left( a_{k,i}\right) _{i\in -%
%TCIMACRO{\U{2115} }%
%BeginExpansion
\mathbb{N}
%EndExpansion
}\in \left\{ 0,1\right\} ^{-%
%TCIMACRO{\U{2115} }%
%BeginExpansion
\mathbb{N}
%EndExpansion
}$ with $a_{k+1}=a_{k}+1$ for each $k\in 
%TCIMACRO{\U{2124} }%
%BeginExpansion
\mathbb{Z}
%EndExpansion
$, such that, for each $k$ with $a_{k}\neq (\ldots ,0,\ldots ,0)$, we turn
right between $X_{k}$\ and $X_{k+1}$\ if $P(a_{k})+Q(a_{k})+R(a_{k})$ is
even, and left otherwise.\ We also consider somme sequences $b_{k}=\left(
b_{k,i}\right) _{i\in -%
%TCIMACRO{\U{2115} }%
%BeginExpansion
\mathbb{N}
%EndExpansion
}$ which satisfy the same properties for $B$.

Each subcurve of $A$ with finitely many segments is contained in some $%
A_{k,m}=\left( X_{k-2^{m}+1},\ldots ,X_{k+2^{m}}\right) $\ with $k\in 
%TCIMACRO{\U{2124} }%
%BeginExpansion
\mathbb{Z}
%EndExpansion
$, $m\in 
%TCIMACRO{\U{2115} }%
%BeginExpansion
\mathbb{N}
%EndExpansion
$\ and $a_{k,i}=0$\ for $i\geq -m+1$.

The sequences which determine the behaviour of $A_{k,m}$ at its turning
points are

\noindent $a_{k-2^{m}}=\left( \ldots ,a_{k-1,i},\ldots ,a_{k-1,-m},0,\ldots
,0\right) $

\noindent ............................................................

\noindent $a_{k-1}=\left( \ldots ,a_{k-1,i},\ldots ,a_{k-1,-m},1,\ldots
,1\right) $

\noindent $a_{k}=\left( \ldots ,a_{k,i},\ldots ,a_{k,-m},0,\ldots ,0\right) $

\noindent ............................................................

\noindent $a_{k+2^{m}-1}=\left( \ldots ,a_{k,i},\ldots ,a_{k,-m},1,\ldots
,1\right) $.

For each $l\in 
%TCIMACRO{\U{2124} }%
%BeginExpansion
\mathbb{Z}
%EndExpansion
$, such that $b_{l,i}=0$\ for $i\geq -m+1$, the sequences which determine
the behaviour of $B_{l,m}=\left( Y_{l-2^{m}+1},\ldots ,Y_{l+2^{m}}\right) $
at its turning points are

\noindent $b_{l-2^{m}}=\left( \ldots ,b_{l-1,i},\ldots ,b_{l-1,-m},0,\ldots
,0\right) $

\noindent ............................................................

\noindent $b_{l-1}=\left( \ldots ,b_{l-1,i},\ldots ,b_{l-1,-m},1,\ldots
,1\right) $

\noindent $b_{l}=\left( \ldots ,b_{l,i},\ldots ,b_{l,-m},0,\ldots ,0\right) $

\noindent ............................................................

\noindent $b_{l+2^{m}-1}=\left( \ldots ,b_{l,i},\ldots ,b_{l,-m},1,\ldots
,1\right) $.

If $b_{l,-m-1}=a_{k,-m-1}$ and $b_{l,-m}=a_{k,-m}$, then $%
b_{l-1,-m-1}=a_{k-1,-m-1}$ and $b_{l-1,-m}=a_{k-1,-m}$. It follows that $%
\left( Y_{l-2^{m}+1},\ldots ,Y_{l}\right) $\ (resp. $\left( Y_{l+1},\ldots
,Y_{l+2^{m}}\right) $) is the image of $\left( X_{k-2^{m}+1},\ldots
,X_{k}\right) $\ (resp. $\left( X_{k+1},\ldots ,X_{k+2^{m}}\right) $)\ by a
translation or a point reflection.

For each $p\in 
%TCIMACRO{\U{2115} }%
%BeginExpansion
\mathbb{N}
%EndExpansion
$, we write $0_{p}=\left( 0,\ldots ,0\right) $\ with $p$ times $0$. We
consider two cases.

If $\left( a_{k,-m-1},a_{k,-m}\right) \neq \left( 0,0\right) $, we consider
the integers $l$ such that

\noindent $b_{l}=\left( \ldots ,b_{l,i},\ldots
,b_{l,-m-8},0,0,0,0,a_{k,-m-3},a_{k,-m-2},a_{k,-m-1},a_{k,-m},0_{m}\right) $.

\noindent For such an $l$, we have

\noindent $b_{l+2^{m+6}+2^{m+5}}=\left( \ldots ,b_{l,i},\ldots
,b_{l,-m-8},0,1,1,0,a_{k,-m-3},a_{k,-m-2},a_{k,-m-1},a_{k,-m},0_{m}\right) $.

\noindent Modulo translations and point reflections, $A_{k,m}$ is equivalent
to $B_{l,m}$ and to $B_{l+2^{m+6}+2^{m+5},m}$. As $%
P(b_{l+2^{m+6}+2^{m+5}})+Q(b_{l+2^{m+6}+2^{m+5}})+R(b_{l+2^{m+6}+2^{m+5}})$
and $P(b_{l})+Q(b_{l})+R(b_{l})$ are not equal, it follows that $A_{k,m}$ is
isomorphic to one of the curves $B_{l,m}$, $B_{l+2^{m+6}+2^{m+5},m}$, and
equivalent to the other one modulo a point reflection.

If $\left( a_{k,-m-1},a_{k,-m}\right) =\left( 0,0\right) $, we consider the
integers $l$ such that

\noindent $b_{l}=\left( \ldots ,b_{l,i},\ldots
,b_{l,-m-8},0,0,0,0,0,1,0_{m+2}\right) $.

\noindent For such an $l$, we have

\noindent $b_{l+2^{m+3}}=\left( \ldots ,b_{l,i},\ldots
,b_{l,-m-8},0,0,0,0,1,1,0_{m+2}\right) $.

\noindent $b_{l+2^{m+6}+2^{m+5}}=\left( \ldots ,b_{l,i},\ldots
,b_{l,-m-8},0,1,1,0,0,1,0_{m+2}\right) $.

\noindent $b_{l+2^{m+6}+2^{m+5}+2^{m+3}}=\left( \ldots ,b_{l,i},\ldots
,b_{l,-m-8},0,1,1,0,1,1,0_{m+2}\right) $.

\noindent Consequently, modulo translations and point reflections, $A_{k,m}$
is equivalent to $B_{l,m}$ and to $B_{l+2^{m+6}+2^{m+5},m}$, or equivalent
to $B_{l+2^{m+3},m}$ and to $B_{l+2^{m+6}+2^{m+5}+2^{m+3},m}$. Moreover, in
each case, $A_{k,m}$ is isomorphic to one of the two curves considered, and
equivalent to the other one modulo a point reflection.

It follows that each curve $D=C_{a,x,g}$\ or $D=C_{g,h}$ is self-avoiding,
since each subcurve of finite length of $D$ is equivalent to a subcurve of
some $C_{n}$, and therefore self-avoiding by Lemma 5.

It remains to be proved that, if a curve $A$ is locally isomorphic to such
curves, then it is isomorphic to one of them. We consider $B=C_{b,y,h}$
which is locally isomorphic to $A$. We write $A=\left( X_{k}\right) _{k\in 
%TCIMACRO{\U{2124} }%
%BeginExpansion
\mathbb{Z}
%EndExpansion
}$ and $B=\left( Y_{k}\right) _{k\in 
%TCIMACRO{\U{2124} }%
%BeginExpansion
\mathbb{Z}
%EndExpansion
}$.

For each $m\in 
%TCIMACRO{\U{2115} }%
%BeginExpansion
\mathbb{N}
%EndExpansion
$, we consider the subcurve $A_{m}=\left( X_{-2^{m}+1},\ldots
,X_{2^{m}}\right) $\ and the subcurves $B_{k,m}=\left( Y_{k-2^{m}+1},\ldots
,Y_{k+2^{m}}\right) $ for $k\in 
%TCIMACRO{\U{2124} }%
%BeginExpansion
\mathbb{Z}
%EndExpansion
$. As $A$ and $B$ are locally isomorphic, each $A_{m}$ is isomorphic to some 
$B_{k_{m},m}$.

By K\"{o}nig's lemma, there exists a subsequence of $\left( b+k_{m}\right)
_{m\in 
%TCIMACRO{\U{2115} }%
%BeginExpansion
\mathbb{N}
%EndExpansion
}$\ which converges to an element $a=\left( a_{i}\right) _{i\in -%
%TCIMACRO{\U{2115} }%
%BeginExpansion
\mathbb{N}
%EndExpansion
}\in \left\{ 0,1\right\} ^{-%
%TCIMACRO{\U{2115} }%
%BeginExpansion
\mathbb{N}
%EndExpansion
}$. If $\left\{ i\in -%
%TCIMACRO{\U{2115} }%
%BeginExpansion
\mathbb{N}
%EndExpansion
\mid a_{i}=0\right\} $ and $\left\{ i\in -%
%TCIMACRO{\U{2115} }%
%BeginExpansion
\mathbb{N}
%EndExpansion
\mid a_{i}=1\right\} $ are infinite, then $A$ is isomorphic to some $%
C_{a,x,g}$. Otherwise, $A$ is isomorphic to some $C_{g_{1},g_{2}}$.~~$%
\blacksquare $\bigskip

\begin{center}
\textbf{References\bigskip }
\end{center}

\noindent \lbrack 1] M. Mend\`{e}s France and J. Shallit, Some planar curves
associated with sums of digits, unpublished.

\noindent \lbrack 2] F. Oger, Paperfolding sequences, paperfolding curves
and local isomorphism, Hiroshima Math. Journal 42 (2012), 37-75.

\noindent \lbrack 3] F. Oger, The number of paperfolding curves in a
covering of the plane, Hiroshima Math. Journal 47 (2017), 1-14.

\noindent \lbrack 4] F. Oger, Coverings of the plane by self-avoiding curves
which satisfy the local isomorphism property, arXiv:2310.19364.\bigskip

\end{document}